\documentclass[11pt]{article}
\usepackage{amssymb}

\usepackage{amsmath}
\usepackage{eufrak}

\newtheorem{theorem}{Theorem}[section]

\newtheorem{corollary}[theorem]{Corollary}

\newtheorem{definition}[theorem]{Definition}
\newtheorem{example}[theorem]{Example}

\newtheorem{lemma}[theorem]{Lemma}

\newtheorem{proposition}[theorem]{Proposition}
\newtheorem{remark}[theorem]{Remark}

\newenvironment{proof}[1][Proof]{\textbf{#1.} }{\ \rule{0.5em}{0.5em}}

\begin{document}

\title{Central limit theorems for multiple Skorohod integrals}
\author{Ivan Nourdin \\
Laboratoire de Probabilit{\'{e}}s et Mod{\`{e}}les Al{\'{e}}atoires\\
Universit{\'{e}} Pierre et Marie Curie\\
Bo{\^{\i}}te courrier 188\\
4 Place Jussieu, 75252 Paris Cedex 5, France\\
\texttt{ivan.nourdin@upmc.fr} \and David Nualart\thanks{%
The work of D. Nualart is supported by the NSF Grant DMS-0604207} \\
Department of Mathematics\\
University of Kansas\\
Lawrence, Kansas 66045, USA\\
\texttt{nualart@math.ku.edu}}
\maketitle

\begin{abstract}
In this paper, we prove a central limit theorem for a sequence of multiple
Skorohod integrals using the techniques of Malliavin calculus. The
convergence is stable, and the limit is a conditionally Gaussian random
variable. Some applications to sequences of multiple stochastic integrals,
and renormalized weighted Hermite variations of the fractional Brownian
motion are discussed. \newline


\end{abstract}

\section{Introduction}

Consider a sequence of random variables $\{F_{n},n\geq 1\}$ defined on a
complete probability space $(\Omega ,\mathcal{F},P)$. Suppose that the 
$\sigma $-field $\mathcal{F}$ is generated by an isonormal Gaussian\ process 
$X=\{X(h),h\in \EuFrak H\}$ on a real separable infinite-dimensional Hilbert space $\EuFrak H$. This
just means that $X$ is a centered Gaussian family of random variables indexed by
the elements of $\EuFrak H$, and such that, for every $h,g\in \EuFrak H$,
\begin{equation}
E\left[ X(h)X(g)\right] =\langle h,g\rangle _{\EuFrak H}.  \label{v1}
\end{equation}
Suppose that the sequence $\{F_{n},n\geq 1\}$ is normalized, that is, 
$E(F_{n})=0$ and $\lim_{n\rightarrow \infty }E(F_{n}^{2})=1$. A natural
problem is to find suitable conditions ensuring that $F_n$ converges in law
towards a given distribution. When the random variables $F_{n}$ belong to
the $q$th Wiener chaos of $X$ (for a fixed $q\geq 2$), then it turns out that the following
conditions are equivalent:

\begin{enumerate}
\item[(i)] $F_n$ converges in law to $N(0,1)$;

\item[(ii)] $\lim_{n\rightarrow \infty }E[F_{n}^{4}]=3$;


\item[(iii)] $\lim_{n\rightarrow \infty }\Vert DF_{n}\Vert _{\EuFrak H}^{2}=q
$ in $L^{2}(\Omega) $.
\end{enumerate}
Here, 
$D$ stands for the
derivative operator in the sense of Malliavin calculus (see Section \ref{Pre} below
for more details).
More precisely, the following bound is in order, where $N$ denotes a standard Gaussian random variable:
\begin{eqnarray}\label{ber-ess}
\sup_{z\in \mathbb{R}}\left| P( F_{n}\leq z) -P(N\leq z)\right|
&\leqslant& \sqrt{E\left[ \left( 1-\frac1q\Vert DF_{n}\Vert _{\EuFrak H}^{2}\right) ^{2}\right] }  \\
\label{4mom}
&\leqslant& \sqrt{\frac{q-1}{3q}}\sqrt{\big|E(F_n^4)-3\big|}.
\end{eqnarray}

The equivalence between conditions (i) and (ii) was proved in Nualart
and Peccati \cite{NP} by means of the Dambis, Dubins and Schwarz
theorem. It implies that the convergence in
distribution of a sequence of multiple stochastic integrals towards a
Gaussian random variable is completely determined by the asymptotic behavior
of their second and fourth moments, which represents a drastic simplification
of the classical ``method of moments and diagrams'' (see, for instance, the
survey by Peccati and Taqqu \cite{Sur}, as well as the references therein).
The equivalence with condition (iii) was proved 
later
by Nualart and
Ortiz-Latorre \cite{NO} using tools of Malliavin calculus. Finally, the
Berry-Esseen's type bound (\ref{ber-ess}) is taken from Nourdin and
Peccati \cite{stein}, while (\ref{4mom}) was shown in Nourdin, Peccati
and Reinert \cite{NPR}. 

\medskip Peccati and Tudor \cite{PT} also obtained a multidimensional
version of the equivalence between (i) and (ii). In particular, they proved that, given
a sequence $\{F_{n},n\geq 1\}$ of $d$-dimensional random vectors such that $%
F_{n}^{i}$ belongs to the $q_{i}$th Wiener chaos for $i=1,\ldots ,d$, where $%
1\leqslant q_{1}\leqslant \ldots \leqslant q_{d}$, then if the covariance
matrix of $F_{n}$ converges to the $d\times d$ identity matrix $I_{d}$, the
convergence in distribution to each component towards the law $N(0,1)$ implies
the convergence in distribution of the whole sequence $F_{n}$ towards the
standard centered Gaussian law $N(0,I_{d})$.

\medskip Recent examples of application of these results are, among others,
the study of $p$-variations of fractional stochastic integrals (Corcuera
\textit{et al.} \cite{CNW}), quadratic functionals of bivariate Gaussian
processes (Deheuvels \textit{et al.} \cite{DPY}), self-intersection local
times of fractional Brownian motion (Hu and Nualart \cite{HN}),
approximation schemes for scalar fractional differential equations
(Neuen\-kirch and Nourdin \cite{NN}), high-frequency CLTs for random fields
on homogeneous spaces (Marinucci and Peccati \cite{MP,MP2} and Peccati \cite%
{Pecc}), needlets analysis on the sphere (Baldi \textit{et al.} \cite{BKMP}%
), estimation of self-similarity orders (Tudor and Viens \cite{TudViens}),
weighted power variations of iterated Brownian motion (Nourdin and
Peccati \cite{IBM}) or bipower variations of Gaussian processes with stationary increments
(Barndorff-Nielsen \textit{et al.} \cite{BN}).

\medskip Since the works by Nualart and Peccati \cite{NP} and Peccati and
Tudor \cite{PT}, great efforts have been made to find similar statements in
the case where the limit is not necessarily Gaussian. In the references \cite%
{PecTaq} and \cite{PecTaqEJP}, Peccati and Taqqu propose sufficient
conditions ensuring that a given sequence of multiple Wiener-It\^{o}
integrals converges stably towards mixtures of Gaussian random
variables. In another direction, Nourdin and Peccati \cite{NouPec} proved an
extension of the above equivalence (i) -- (iii) for a sequence of random
variables $\{F_{n},n\geq 1\}$ in a fixed $q$th Wiener chaos, $q\geq 2$,
where the limit law is $2\,G_{\nu /2}-\nu $, $G_{\nu /2}$ being the Gamma
distribution with parameter $\nu /2$. 

\medskip The purpose of the present paper is to study the convergence in
distribution of a sequence of random variables of the form $F_{n}=\delta
^{q}(u_{n})$, where $u_{n}$ are random variables with values in $\EuFrak %
H^{\otimes q}$ (the $q$th tensor product of $\EuFrak H$) and $\delta ^{q}$ denotes the 
multiple Skorohod integral (that is, $\delta^2(u)=\delta(\delta(u))$,
$\delta^3(u)=\delta(\delta(\delta(u)))$, and so on), towards a mixture of Gaussian random
variables. Our main abstract result, Theorem \ref{main-thm}, roughly says that under some
technical
conditions, if  $\left\langle u_{n},D^{q}F_{n}\right\rangle _{\EuFrak %
H^{\otimes q}}$ converges in $L^{1}(\Omega)$ to a nonnegative random variable $S^{2}$%
, then the sequence $F_{n}$ converges stably to a random variable $F$ with
conditional characteristic function $E\left( e^{i\lambda F}\right|
X)=E\big( e^{-\frac{\lambda ^{2}}{2}S^{2}}\big) $. Notice that if $u_{n}$
is deterministic, then $F_{n}$ belongs to the $q$th Wiener chaos, and we
have a sequence of the type considered above. In particular, if $S^{2}$ is also
deterministic, we recover the fact that condition (iii) above implies the
convergence in distribution to the law $N(0,1)$.

\medskip We develop some particular applications of Theorem \ref{main-thm}
in the following directions. First, we consider a sequence of random
variables in a fixed Wiener chaos and we derive new criteria for the
convergence to a mixture of Gaussian laws. Second, we show the convergence
in law of the sequence $\delta ^{q}(u_{n})$, where $q\ge 2$ and  $u_{n}$ is a  
$q$-parameter process of the form
\begin{equation*}
u_{n}=n^{qH-\frac{1}{2}}\sum_{k=0}^{n-1}f(B_{k/n})\mathbf{1}%
_{(k/n,(k+1)/n]^{q}},
\end{equation*}%
towards the random variable $\sigma _{H,q}\int_{0}^{1}f(B_{s})dW_{s}$, where $B$ is
a fractional Brownian motion with Hurst parameter $H\in \left( \frac{1}{4q},%
\frac{1}{2}\right) $, $W$ is a standard Brownian motion independent of $B$, and $\sigma_{H,q}$ denotes some
positive constant. This convergence allows us to establish a new asymptotic result for
the behavior of the weighted $q$th Hermite variation of the fractional Brownian
motion with Hurst parameter $H\in \left( \frac{1}{4q},\frac{1}{2}\right) $,
which complements and provides a new perspective to the results proved
by Nourdin \cite{No1}, 
Nourdin, Nualart and Tudor \cite{NNT}, and
Nourdin and R\'{e}veillac \cite{NR}. The reader is referred to Section \ref{fbm}
for a detailed description of these results.

\medskip The paper is organized as follows. In Section \ref{Pre}, we present
some preliminary results about Malliavin calculus. 
Section \ref{main} contains the statement and the proof of the main abstract result.
In Section \ref{sect-multiple}, we apply it to sequences
of multiple stochastic integrals, while Section \ref{fbm} 
focuses on the applications to the weighted Hermite variations of the
fractional Brownian motion.

\setcounter{equation}{0}

\section{Preliminaries}

\label{Pre}

Let $\EuFrak H$ be a real separable infinite-dimensional Hilbert space. For any integer $q\geq 1$, let $%
\EuFrak H^{\otimes q}$ be the $q$th tensor product of $\EuFrak H$. Also, we denote
by $\EuFrak H^{\odot q}$ the $q$th symmetric tensor product.

Suppose that $X=\{X(h),\,h\in \EuFrak H\}$ is an isonormal Gaussian process on
$\EuFrak H$, defined on some probability space $(\Omega ,\mathcal{F},P)$.
Recall that this means that the covariance of $X$ is given in terms of the
scalar product of $\EuFrak H$ by (\ref{v1}). Assume from now on that $\mathcal{F}$ is
generated by $X$.

For every integer $q\geq 1$, let $\mathcal{H}_{q}$ be the $q$th Wiener chaos of $X$,
that is, the closed linear subspace of $L^{2}(\Omega)$
generated by the random variables $\{H_{q}(X(h)),h\in \EuFrak H,\left\|
h\right\| _{\EuFrak H}=1\}$, where $H_{q}$ is the $q$th Hermite polynomial defined by
$$
H_q(x)=\frac{(-1)^q}{q!}e^{x^2/2}\frac{d^q}{dx^q}\big(e^{-x^2/2}\big).
$$
We denote by $\mathcal{H}_{0}$ the space of constant random variables. For
any $q\geq 1$, the mapping $I_{q}(h^{\otimes q})=q!H_{q}(X(h))$ provides a
linear isometry between $\EuFrak H^{\odot q}$
(equipped with the modified norm $\sqrt{q!}\left\| \cdot \right\| _{\EuFrak %
H^{\otimes q}}$) and $\mathcal{H}_{q}$ (equipped with the $L^2(\Omega)$ norm). For $q=0$, by convention $%
\mathcal{H}_{0}=\mathbb{R}$, and $I_{0}$ is the identity map.

It is well-known (Wiener chaos expansion) that $L^{2}(\Omega)$
can be decomposed into the infinite orthogonal sum of the spaces $\mathcal{H}%
_{q}$. That is, any square integrable random variable $F\in L^{2}(\Omega)$ admits the following chaotic expansion:
\begin{equation}
F=\sum_{q=0}^{\infty }I_{q}(f_{q}),  \label{E}
\end{equation}%
where $f_{0}=E[F]$, and the $f_{q}\in \EuFrak H^{\odot q}$, $q\geq 1$, are
uniquely determined by $F$. For every $q\geq 0$, we denote by $J_{q}$ the
orthogonal projection operator on the $q$th Wiener chaos. In particular, if $%
F\in L^{2}(\Omega)$ is as in (\ref{E}), then $%
J_{q}F=I_{q}(f_{q})$ for every $q\geq 0$.

Let $\{e_{k},\,k\geq 1\}$ be a complete orthonormal system in $\EuFrak H$.
Given $f\in \EuFrak H^{\odot p}$, $g\in \EuFrak H^{\odot q}$ and $%
r\in\{0,\ldots ,p\wedge q\}$, the $r$th contraction of $f$ and $g$ 
is the element of $\EuFrak H^{\otimes (p+q-2r)}$ defined by
\begin{equation}
f\otimes _{r}g=\sum_{i_{1},\ldots ,i_{r}=1}^{\infty }\langle
f,e_{i_{1}}\otimes \ldots \otimes e_{i_{r}}\rangle _{\EuFrak H^{\otimes
r}}\otimes \langle g,e_{i_{1}}\otimes \ldots \otimes e_{i_{r}}\rangle _{%
\EuFrak H^{\otimes r}}.  \label{v2}
\end{equation}%
Notice that $f\otimes _{r}g$ is not necessarily symmetric. We denote its
symmetrization by $f\widetilde{\otimes }_{r}g\in \EuFrak H^{\odot (p+q-2r)} $%
. Moreover, $f\otimes _{0}g=f\otimes g$ equals the tensor product of $f$ and
$g$ while, for $p=q$, $f\otimes _{q}g=\langle f,g\rangle _{\EuFrak %
H^{\otimes q}}$.

In the particular case $\EuFrak H=L^{2}(A,\mathcal{A},\mu )$, where $%
(A,\mathcal{A})$ is a measurable space and $\mu $ is a $\sigma $-finite and
non-atomic measure, one has that $\EuFrak H^{\odot q}=L_{s}^{2}(A^{q},%
\mathcal{A}^{\otimes q},\mu ^{\otimes q})$ is the space of symmetric and
square integrable functions on $A^{q}$. Moreover, for every $f\in \EuFrak %
H^{\odot q}$, $I_{q}(f)$ coincides with the multiple Wiener-It\^{o} integral
of order $q$ of $f$ with respect to $X$ (introduced by It\^{o} in \cite{It})
and (\ref{v2}) can be written as
\begin{eqnarray*}
&&(f\otimes _{r}g)(t_1,\ldots,t_{p+q-2r}) 
=\int_{A^{r}}f(t_{1},\ldots ,t_{p-r},s_{1},\ldots ,s_{r}) \\
&&\hskip2cm \times \ g(t_{p-r+1},\ldots ,t_{p+q-2r},s_{1},\ldots ,s_{r})d\mu
(s_{1})\ldots d\mu (s_{r}).
\end{eqnarray*}

Let us now introduce some basic elements of the Malliavin calculus with respect
to the isonormal Gaussian process $X$. We refer the reader to Nualart \cite%
{Nbook} for a more detailed presentation of these notions. Let $\mathcal{S}$
be the set of all smooth and cylindrical random variables of
the form
\begin{equation}
F=g\left( X(\phi _{1}),\ldots ,X(\phi _{n})\right) ,  \label{v3}
\end{equation}%
where $n\geq 1$, $g:\mathbb{R}^{n}\rightarrow \mathbb{R}$ is a infinitely
differentiable function with compact support, and $\phi _{i}\in \EuFrak H$.
The Malliavin derivative of $F$ with respect to $X$ is the element of $%
L^{2}(\Omega ,\EuFrak H)$ defined as
\begin{equation*}
DF\;=\;\sum_{i=1}^{n}\frac{\partial g}{\partial x_{i}}\left( X(\phi
_{1}),\ldots ,X(\phi _{n})\right) \phi _{i}.
\end{equation*}
By iteration, one can
define the $q$th derivative $D^{q}F$ for every $q\geq 2$, which is an element of $L^{2}(\Omega ,%
\EuFrak H^{\odot q})$.

For $q\geq 1$ and $p\geq 1$, ${\mathbb{D}}^{q,p}$ denotes the closure of $%
\mathcal{S}$ with respect to the norm $\Vert \cdot \Vert_{\mathbb{D}^{q,p}}$, defined by
the relation
\begin{equation*}
\Vert F\Vert _{\mathbb{D}^{q,p}}^{p}\;=\;E\left[ |F|^{p}\right] +\sum_{i=1}^{q}E\left(
\Vert D^{i}F\Vert _{\EuFrak H^{\otimes i}}^{p}\right) .
\end{equation*}
The Malliavin derivative $D$ verifies the following chain rule. If $%
\varphi :\mathbb{R}^{n}\rightarrow \mathbb{R}$ is continuously
differentiable with bounded partial derivatives and if $F=(F_{1},\ldots
,F_{n})$ is a vector of elements of ${\mathbb{D}}^{1,2}$, then $\varphi
(F)\in {\mathbb{D}}^{1,2}$ and
\begin{equation*}
D\varphi (F)=\sum_{i=1}^{n}\frac{\partial \varphi }{\partial x_{i}}%
(F)DF_{i}.
\end{equation*}

We denote by $\delta $ the adjoint of the operator $D$, also called the 
divergence operator. The operator $\delta $ is also called the
Skorohod integral because in the case of the Brownian motion it coincides
with the anticipating stochastic integral introduced by Skorohod in \cite{Sk}. 
A random element $u\in L^{2}(\Omega ,\EuFrak %
H)$ belongs to the domain of $\delta $, noted $\mathrm{Dom}\delta $, if and
only if it verifies
\begin{equation*}
\big|E\big(\langle DF,u\rangle _{\EuFrak H}\big)\big|\leq c_{u}\,\sqrt{E(F^2)}
\end{equation*}%
for any $F\in \mathbb{D}^{1,2}$, where $c_{u}$ is a constant depending only
on $u$. If $u\in \mathrm{Dom}\delta $, then the random variable $\delta (u)$
is defined by the duality relationship (called `integration by parts
formula'):
\begin{equation}
E(F\delta (u))=E\big(\langle DF,u\rangle _{\EuFrak H}\big),  \label{ipp}
\end{equation}%
which holds for every $F\in {\mathbb{D}}^{1,2}$. The formula (\ref%
{ipp}) extends to the multiple Skorohod integral $\delta ^{q}$, and we
have
\begin{equation}
E\left( F\delta ^{q}(u)\right) =E\left( \left\langle D^{q}F,u\right\rangle _{%
\EuFrak H^{\otimes q}}\right)  \label{dual}
\end{equation}%
for any element $u$ in the domain of $\delta ^{q}$ and any random variable $%
F\in \mathbb{D}^{q,2}$. Moreover, $\delta ^{q}(h)=I_{q}(h)$ for any $h\in %
\EuFrak H^{\odot q}$. 

The following property will be extensively used in the paper.

 \begin{lemma}
Let $q\geq 1$ be an integer. Suppose that  $F\in {\mathbb{D}}^{q,2}$, and
let $u$ be a symmetric element in $\mathrm{Dom}\delta ^{q}$. Assume that,
for any $\ 0\leq r+j\leq q$,
$\left\langle D^{r}F,\delta ^{j}(u)\right\rangle _{\EuFrak H^{\otimes r}}\in
L^{2}(\Omega ,\EuFrak H^{\otimes q-r-j})$. Then, for any  $r=0,\ldots ,q-1$, $%
\left\langle D^{r}F,u\right\rangle _{\EuFrak H^{\otimes r}}$ belongs to the
domain of $\delta ^{q-r}$ and we have 
\begin{equation}
F\delta ^{q}(u)=\sum_{r=0}^{q}\binom{q}{r}\delta ^{q-r}\left( \left\langle
D^{r}F,u\right\rangle _{\EuFrak H^{\otimes r}}\right) .  \label{t3}
\end{equation}
(We use the convention that $\delta^0(v)=v$, $v\in\mathbb{R}$, and $D^0F=F$, $F\in L^2(\Omega)$.)
\end{lemma}

\begin{proof}
We prove this lemma by induction on $q$. For $q=1$ it 
reads $F\delta(u)=\delta(Fu)+\langle DF,u\rangle_\EuFrak H$, and this formula is well-known,
see e.g.
\cite[Proposition 1.3.3]{Nbook}. Suppose the result is true for $q$. Then,
if $u$ belongs to the domain of $\delta ^{q+1}$,  by the induction
hypothesis applied to $\delta (u)$,%
\begin{equation}
F\delta ^{q+1}(u)=F\delta ^{q}(\delta (u))=\sum_{r=0}^{q}\binom{q}{r}\delta
^{q-r}\left( \left\langle D^{r}F,\delta (u)\right\rangle _{\EuFrak %
H^{\otimes r}}\right) .  \label{t1}
\end{equation}%
On the other hand%
\begin{equation}
\left\langle D^{r}F,\delta (u)\right\rangle _{\EuFrak H^{\otimes r}}=\delta
\left( \left\langle D^{r}F, u \right\rangle _{\EuFrak H^{\otimes
r}}\right) +\left\langle D^{r+1}F,u\right\rangle _{\EuFrak H^{\otimes r}}.
\label{t2}
\end{equation}%
Finally, substituting (\ref{t2}) into (\ref{t1}) yields the desired result.
\end{proof}
\medskip

For any Hilbert space $V$, we denote by $\mathbb{D}^{k,p}(V)$ the
corresponding Sobolev space of $V$-valued random variables (see \cite[page 31]{Nbook}%
). The operator $\delta^q$   is continuous from $\mathbb{D}^{k,p}(\EuFrak H^{\otimes q})$
to $\mathbb{D}^{k-q,p}$, for any $p>1$ and  any integers $k\ge q\ge 1$, that is,
we have%
\begin{equation}
\left\| \delta^q (u)\right\| _{\mathbb{D}^{k-q,p}}\leq c_{k,p}\left\| u\right\| _{\mathbb{D}%
^{k,p}(\EuFrak H^{\otimes q})}  \label{Me2}
\end{equation}%
for all $u\in\mathbb{D}^{k,p}(\EuFrak H^{\otimes q})$, and some constant $c_{k,p}>0$.
These estimates are consequences of Meyer inequalities (see %
\cite[Proposition 1.5.7]{Nbook}).
 In particular, these estimates imply that  
$\mathbb{D}^{q,2}(\EuFrak H^{\otimes q})\subset \mathrm{Dom}\delta ^{q}$ for any integer $q\ge 1$.

We will also use the following commutation relationship between the
Malliavin derivative and the Skorohod integral (see \cite[Proposition 1.3.2]{Nbook})%
\begin{equation}
D\delta (u)=u+\delta (Du),  \label{comm1}
\end{equation}%
for any $u\in \mathbb{D}^{2,2}(\EuFrak H)$. By induction we can show
the following formula for any symmetric element $u$ in $\mathbb{D}^{j+k,2}(%
\EuFrak H^{\otimes j})$ 
\begin{equation}
D^{k}\delta ^{j}(u)=\sum_{i=0}^{j\wedge k}\binom{k}{i}\binom{j}{i}i!\delta
^{j-i}(D^{k-i}u).  \label{t5}
\end{equation}%
We will make use of the following formula for the variance of a multiple
Skorohod integral. Let $u,v\in\mathbb{D}^{2q,2}(\EuFrak H^{\otimes q})\subset{\rm Dom}\delta^q$
be two symmetric functions. Then
\begin{eqnarray}
E(\delta ^{q}(u)\delta ^{q}(v)) &=&E(\left\langle u,D^{q}(\delta
^{q}(v))\right\rangle _{\EuFrak H^{\otimes q}})  \notag \\
&=&\sum_{i=0}^{q}\binom{q}{i}^{2}i!E\left( \left\langle u,\delta
^{q-i}(D^{q-i}v)\right\rangle _{\EuFrak H^{\otimes q}}\right)  \notag \\
&=&\sum_{i=0}^{q}\binom{q}{i}^{2}i!E\left(\left\langle D^{q-i}u,D^{q-i}v\right\rangle_{\EuFrak H^{\otimes (2q-i)}}\right) .  \label{t13}
\end{eqnarray}

The operator $L$ is defined on the Wiener chaos expansion as
\begin{equation*}
L=\sum_{q=0}^{\infty }-qJ_{q},
\end{equation*}
and is called the infinitesimal generator of the Ornstein-Uhlenbeck
semigroup. The domain of this operator in $L^{2}(\Omega)$ is the set%
\begin{equation*}
\mathrm{Dom}L=\{F\in L^{2}(\Omega ):\sum_{q=1}^{\infty }q^{2}\left\|
J_{q}F\right\| _{L^{2}(\Omega )}^{2}<\infty \}=\mathbb{D}^{2,2}\text{.}
\end{equation*}%
There is an important relation between the operators $D$, $\delta $ and $L$
(see \cite[Proposition 1.4.3]{Nbook}). A random variable $F$ belongs to the
domain of $L$ if and only if $F\in \mathrm{Dom}\left( \delta D\right) $
(i.e. $F\in {\mathbb{D}}^{1,2}$ and $DF\in \mathrm{Dom}\delta $), and in
this case
\begin{equation}
\delta DF=-LF.  \label{k1}
\end{equation}

Note also that a random variable $F$ as in (\ref{E}) is in ${\mathbb{D}}%
^{1,2}$ if and only if
\begin{equation*}
\sum_{q=1}^{\infty }qq!\Vert f_{q}\Vert _{\EuFrak H^{\otimes q}}^{2}<\infty ,
\end{equation*}%
and, in this case, $E\left( \Vert DF\Vert _{\EuFrak H}^{2}\right)
=\sum_{q\geq 1}qq!\Vert f_{q}\Vert _{\EuFrak H^{\otimes q}}^{2}$. If $\EuFrak H=%
L^{2}(A,\mathcal{A},\mu )$ (with $\mu $ non-atomic), then the
derivative of a random variable $F$ as in (\ref{E}) can be identified with
the element of $L^{2}(A\times \Omega )$ given by
\begin{equation}
D_{a}F=\sum_{q=1}^{\infty }qI_{q-1}\left( f_{q}(\cdot ,a)\right) ,\quad a\in
A.  \label{dtf}
\end{equation}

Finally, we need the definition of stable convergence (see, for instance, the
original paper \cite{reniy}, or the book \cite{JacSh} for an exhaustive
discussion of stable convergence).

\begin{definition}
\label{def1} 
Let $F_{n}$ be a sequence of random variables defined on the probability space
$(\Omega,\mathcal{F},P)$, and suppose that $F$ is a random variable defined on an enlarged
probability space $(\Omega ,\mathcal{G},P)$, with $\mathcal{F}\subseteq
\mathcal{G}$. 
We say that $F_{n}$ converges
$\mathcal{G}${\rm -stably} to $F$ (or only {\rm stably} when the context is clear)
if,
for any continuous and bounded function $f:\mathbb{R}%
\rightarrow \mathbb{R}$ and any bounded $\mathcal{F}$-measurable random
variable $Z$, we have
$
E\left[ f(F_{n})Z\right] \rightarrow E\left[ f(F)Z\right]
$
as $n$ tends to infinity.
\end{definition}

\setcounter{equation}{0}

\section{Convergence in law of multiple Skorohod integrals}

As in the previous section, $X=\{X(h), h\in \EuFrak H\}$ is an isonormal
Gaussian process associated with a real separable infinite-dimensional Hilbert space $\EuFrak H$. %
\label{main} The next theorem is the main abstract result of the present paper.

\begin{theorem}
\label{main-thm} Fix an integer $q\geq 1$, and suppose that $F_{n}$ is a sequence
of random variables of the form $F_{n}=\delta ^{q}(u_{n})$, 
for some symmetric functions $u_n $ in  
$\mathbb{D}^{2q,2q}(\EuFrak {H}^{\otimes q})$.
Suppose
moreover
that the sequence $F_{n}$ is bounded in $L^{1}(\Omega)$, and that:

\begin{enumerate}
\item[(i)] $\left\langle u_{n},(DF_{n})^{\otimes k_{1}}\otimes 
\ldots \otimes (D^{q-1}F_{n})^{\otimes k_{q-1}}
\otimes h 
\right\rangle _{\EuFrak H^{\otimes q}}$ converges in $L^{1}(\Omega)$ to zero,
for all integers $r,k_{1},\ldots ,k_{q-1}\geq 0$ such that 
$$
k_{1}+2k_{2}+\ldots +(q-1)k_{q-1}+r=q,
$$ 
and all $h\in \EuFrak H^{\otimes r}$;

\item[(ii)] $\left\langle u_{n},D^{q}F_{n} \right\rangle _{\EuFrak %
H^{\otimes q}}$ converges in $L^{1}(\Omega)$ to a nonnegative random variable $S^{2}$%
.
\end{enumerate}

Then, $F_{n}$ converges stably to a random variable with conditional
Gaussian law $N(0,S^2)$ given $X$.
\end{theorem}

\begin{remark}\label{rk31} 
{\rm 
When $q=1$, condition (i) of the theorem is that $%
\langle u_{n},h\rangle _{\EuFrak H}$ converges to zero in $L^{1}(\Omega)$, for each $%
h\in \EuFrak H$. When $q=2$, condition (i) means that $\langle
u_{n},h\otimes g\rangle _{\EuFrak H^{\otimes 2}}$, $\langle u_{n},
DF_{n}\otimes h\rangle _{\EuFrak H^{\otimes 2}}$ and $\langle
u_{n},DF_{n}\otimes DF_{n}\rangle _{\EuFrak H^{\otimes 2}}$ converge to zero
in $L^{1}(\Omega)$, for each $h,g\in \EuFrak H$. And so on. 
} 
\end{remark}

\noindent
\begin{proof}[Proof of Theorem \ref{main-thm}]
Taking into account Definition \ref{def1}, it suffices to show that for any $%
h_{1},\ldots ,h_{m}\in \EuFrak H$, the sequence $$\xi
_{n}=(F_{n},X(h_{1}),\ldots ,X(h_{m}))$$ converges in distribution to
a vector
$
(F_{\infty },X(h_{1}),\ldots ,X(h_{m})),
$
where $F_{\infty }$ satisfies, for any $\lambda \in \mathbb{R}$,
\begin{equation}
E(e^{i\lambda F_{\infty }}|X(h_{1}),\ldots ,X(h_{m}))=e^{-\frac{\lambda ^{2}%
}{2}S^{2}}.  \label{bb1}
\end{equation}%
Since the sequence $F_{n}$ is bounded in $%
L^{1}(\Omega)$, the sequence $\xi _{n}$ is tight. Assume that $(F_{\infty
},X(h_{1}),\ldots ,X(h_{m}))$ denotes the limit in law of a certain
subsequence of $\xi _{n}$, denoted again by $\xi _{n}$.

Let $Y=\phi (X(h_{1}),\ldots ,X(h_{m}))$, with $\phi \in \mathcal{C}%
_{b}^{\infty }(\mathbb{R}^{m})$ ($\phi $ is infinitely differentiable,
bounded, with bounded partial derivatives of all orders), and consider $\phi
_{n}(\lambda )=E\left( e^{i\lambda F_{n}}Y\right) $ for $\lambda \in \mathbb{%
R}$. The convergence in law of $\xi _{n}$, together with the fact that $F_{n}$ is
bounded in $L^{1}(\Omega)$, imply that
\begin{equation}
\lim_{n\rightarrow \infty }\phi _{n}^{\prime }(\lambda )=\lim_{n\rightarrow
\infty }iE\left( F_{n}e^{i\lambda F_{n}}Y\right) =\ iE(F_{\infty
}e^{i\lambda F_{\infty }}Y).   \label{a134}
\end{equation}%
On the other hand, by (\ref{dual}) and the Leibnitz rule for $D^q$, we obtain%
\begin{eqnarray*}
\phi _{n}^{\prime }(\lambda ) &=&iE(F_{n}e^{i\lambda F_{n}}Y)=iE\left(
\delta ^{q}(u_{n})e^{i\lambda F_{n}}Y\right) \\
&=&iE\left( \left\langle
u_{n},D^{q}\left( e^{i\lambda F_{n}}Y\right) \right\rangle _{\EuFrak %
H^{\otimes q}}\right) \\
&=&i\sum_{a=0}^{q}\binom{q}{a}E\left( \left\langle u_{n},D^{a}\left(
e^{i\lambda F_{n}}\right) \widetilde{\otimes }D^{q-a}Y\right\rangle _{%
\EuFrak H^{\otimes q}}\right) \\
&=&i\sum_{a=0}^{q}\binom{q}{a}\sum \frac{a!}{k_{1}!\ldots k_{a}!}(i\lambda
)^{k_{1}+\cdots +k_{a}} \\
&&\times E\left( e^{i\lambda F_{n}}\left\langle u_{n},(DF_{n})^{\otimes
k_{1}}\widetilde{\otimes }\ldots \widetilde{\otimes }(D^{a}F_{n})^{\otimes
k_{a}}\widetilde{\otimes }D^{q-a}Y\right\rangle _{\EuFrak H^{\otimes
q}}\right) \\
&=&i\sum_{a=0}^{q}\binom{q}{a}\sum \frac{a!}{k_{1}!\ldots k_{a}!}(i\lambda
)^{k_{1}+\cdots +k_{a}} \\
&&\times E\left( e^{i\lambda F_{n}}\left\langle u_{n},(DF_{n})^{\otimes
k_{1}}\otimes \ldots \otimes (D^{a}F_{n})^{\otimes
k_{a}}\otimes D^{q-a}Y\right\rangle _{\EuFrak H^{\otimes
q}}\right) ,
\end{eqnarray*}%
where the second sum in the two last equalities runs over all sequences of
integers $(k_{1},\ldots ,k_{a})$ such that $k_{1}+2k_{2}+\ldots +ak_{a}=a$,
due to the Fa\'{a} di Bruno's formula. \ By condition (i), this yields that
\begin{equation*}
\phi _{n}^{\prime }(\lambda )=-\lambda E\left( e^{i\lambda
F_{n}}\left\langle u_{n},D^{q}F_{n}\right\rangle _{\EuFrak H^{\otimes
q}}Y\right) +R_{n},
\end{equation*}%
with $R_{n}$ converging to zero as $n\rightarrow \infty $. Using condition
(ii) and (\ref{a134}), we obtain that
\begin{equation*}
iE(F_{\infty }e^{i\lambda F_{\infty }}Y)=-\lambda E\left( e^{i\lambda
F_{\infty }}S^{2}Y\right) .
\end{equation*}%
Since $S^{2}$ is defined through condition (ii), it is in particular
measurable with respect to $X$. Thus, the following linear differential
equation verified by the conditional characteristic function of $F_{\infty }$
holds:
\begin{equation*}
\frac{\partial }{\partial \lambda }E(e^{i\lambda F_{\infty
}}|X(h_{1}),\ldots ,X(h_{m}))=-\lambda \,S^{2}\,E(e^{i\lambda F_{\infty
}}|X(h_{1}),\ldots ,X(h_{m})).
\end{equation*}%
By solving it, we obtain (\ref{bb1}), which yields the desired conclusion.
\end{proof}

\medskip
 
The next corollary provides stronger but easier conditions for the stable
convergence.

\begin{corollary}
\label{cor-1} For a fixed $q\geq 1$, suppose that $F_{n}$ is a sequence of
random variables of the form $F_{n}=\delta ^{q}(u_{n})$, 
for some symmetric functions $u_n $ in  
$\mathbb{D}^{2q,2q}(\EuFrak {H}^{\otimes q})$. 
Suppose moreover
that the sequence $F_{n}$ is bounded in $\mathbb{D}^{q,p}$ for all $p\geq 2$%
, and that:

\begin{enumerate}
\item[(i')] $\langle u_{n},h\rangle _{\EuFrak H^{\otimes q}}$ converges to
zero in $L^{1}(\Omega)$ for all $h\in \EuFrak H^{\otimes q}$; and $u_{n}\otimes
_{l }D^{l }F_{n}$ converges to zero in $L^{2}(\Omega ;\EuFrak %
H^{\otimes (q-l )})$ for all $l =1,\ldots ,q-1$;

\item[(ii)] $\left\langle u_{n},D^{q}F_{n} \right\rangle _{\EuFrak %
H^{\otimes q}}$ converges in $L^{1}(\Omega)$ to a nonnegative random variable $S^{2}$%
.
\end{enumerate}

Then, $F_{n}$ converges stably to a random variable with conditional
Gaussian law $N(0,S^2)$ given $X$.
\end{corollary}

\begin{proof}
\noindent \quad It suffices to show that condition (i') implies condition
(i) in Theorem \ref{main-thm}. When $k_{a}\not=0$ 
for $1\leq a\leq q-1$, we have, for all $h\in\EuFrak H^{\otimes r}$ (with $r=q-k_1-2k_2-\ldots-ak_a$),
\begin{eqnarray*}
&&\left| \left\langle u_{n},(DF_{n})^{\otimes k_{1}}\otimes
\ldots \otimes (D^{a}F_{n})^{\otimes k_{a}}\otimes 
h\right\rangle _{\EuFrak H^{\otimes q}}\right| \\
&=&\bigg|\bigg\langle u_{n}\otimes_a D^aF_n, \\
&&\hskip.5cm\left.\left.(DF_{n})^{\otimes k_{1}}\otimes
\ldots\otimes (D^{a-1}F_{n})^{\otimes k_{a-1}}\otimes
 (D^{a}F_{n})^{\otimes (k_{a}-1)}\otimes
h\right\rangle _{\EuFrak H^{\otimes (q-a)}}\right| \\
&\leq & \left\| u_{n}\otimes _{a}D^{a}F_{n}\right\| _{\EuFrak H^{\otimes
(q-a)}}  \\
&&\hskip.1cm\times \left\| (DF_{n})^{\otimes k_{1}}\otimes \ldots \otimes
(D^{a-1}F_{n})^{\otimes k_{a-1}}\otimes
(D^{a}F_{n})^{\otimes (k_{a}-1)}\otimes h\right\| _{\EuFrak H^{\otimes
(q-a)}}.
\end{eqnarray*}%
The second factor is bounded in $L^{2}(\Omega)$, and the first factor converges to
zero in $L^{2}(\Omega)$, for all $a=1,\dots ,q-1$. In the case $a=0$ we have that $%
\langle u_{n},h\rangle _{\EuFrak H^{\otimes q}}$ converges to zero in $L^{1}(\Omega)$%
, for all $h\in \EuFrak H^{\otimes q}$, by condition (i'). This completes
the proof.
\end{proof}

\section{Multiple stochastic integrals}\label{sect-multiple}
Suppose that $\EuFrak H$ is a Hilbert space $L^2(A, \mathcal{A}, \mu)$,
where $%
(A,\mathcal{A})$ is a measurable space and $\mu $ is a $\sigma $-finite and
non-atomic measure.

Fix an integer $m\geq 2$, and consider a sequence of multiple stochastic
integrals $\{F_{n}=I_{m}(g_{n}),\,n\geq 1\}$ with $g_{n}\in \EuFrak H^{\odot
m}$.   
We would like to apply Theorem \ref{main-thm}  with $q=1$ to the sequence $F_n$. 
To do this, we represent each $F_n$ as
$$F_n=\delta(u_n),\quad\mbox{with $u_n=I_{m-1} (\widehat{g}_n)$},$$ 
for 
$\widehat{g}_n\in  \EuFrak  H^{\otimes m}$ some function 
which is symmetric in the first $m-1$ variables.  

Notice that, from (\ref{dtf}), we have $DF_{n}=mI_{m-1}(g_{n})$. 
Hence, since $F_{n}=-\frac{1}{m}LF_{n}=\frac{1}{m}\delta (DF_{n})$
by (\ref{k1}), $g_n$ is always 
a possible choice for 
$\widehat{g}_n$. 
(In this case, $\widehat{g}_n$ is symmetric in all the variables.)
However, as observed, for instance, in Example \ref{expoire} below, 
the choice $\widehat{g}_n =g_n $ does not allow to conclude in general.

\begin{proposition}
\label{prop-2} For a fixed integer $m\geq 2$, let $F_{n}$ be a sequence of
random variables of the form $F_{n}=I_{m}(g_{n})$, with $g_{n}\in \EuFrak %
H^{\odot m}$. Suppose moreover that $F_n$ is bounded in $L^2(\Omega)$ and that 
$F_n=\delta(u_n)$, where $u_n=I_{m-1} (\widehat{g}_n)$, for $\widehat{g}_n\in  \EuFrak  H^{\otimes m}$ 
some function which is 
symmetric  in the first $m-1$ variables. Finally, assume that:

\begin{enumerate}
\item[(a)] $\langle \widehat{g}_n\otimes _{m-1}\widehat{g}_n,h^{\otimes 2}\rangle _{\EuFrak %
H^{\otimes 2}}$ converges to zero for all $h\in \EuFrak H$;

\item[(b)] $ \langle u_n, DF_{n} \rangle_{\EuFrak H} $ converges in $%
L^{1}(\Omega)$ to a non negative random variable $S^{2}$.
\end{enumerate}

Then, $F_{n}$ converges stably to a random variable with conditional
Gaussian law $N(0,S^2)$ given $X$.
\end{proposition}

\begin{proof}
\noindent It suffices to apply Theorem \ref{main-thm} to $%
u_{n}=I_{m-1}(\widehat{g}_n)$ and $q=1$. 
Indeed, we have
\begin{eqnarray*}
E\left( \langle u_{n},h\rangle _{\EuFrak H}^{2}\right) &=&E\left(
\left\langle I_{m-1}(\widehat{g}_n),h\right\rangle _{\EuFrak H}^{2}\right) =E\left(
I_{m-1}(\widehat{g}_n\otimes _{1}h)^{2}\right) \\
&=&(m-1)!\left\| \widehat{g}_n\otimes _{1}h\right\| _{\EuFrak H^{\otimes (m-1)}}^{2}\\
&=&(m-1)!\left\langle \widehat{g}_n\otimes _{m-1}\widehat{g}_n,h^{\otimes 2}\right\rangle _{%
\EuFrak H^{\otimes 2}}\rightarrow 0,
\end{eqnarray*}%
which implies condition (i) in Theorem \ref{main-thm}, see also Remark \ref{rk31}%
. Condition (ii) in Theorem \ref{main-thm} follows from (b).
\end{proof}

\begin{example}\label{expoire}
{\rm  
(see also \cite[Proposition 2.1]{PecYor} or \cite[Proposition 18]{PecTaq} for two different proofs using other techniques).
Suppose that $\{W_t,\ t\in [0,1]\}$ is a standard Brownian motion. 
(This corresponds to $A=[0,1]$ and $\mu$ the Lebesgue measure.)  
Assume that $m=2$ and take $g_n(s,t)=\frac 12 \sqrt{n} (s\vee t)^n$.  Then
\[
F_n = I_2(g_n)= \sqrt{n} \int_0^1  t^n W_t dW_t,
\]
and
\[
D_sF_n = \sqrt{n} s^n W_s + \sqrt{n} \int_s^1 t^n W_t dW_t.
\]
We can take $u_n(t)=\sqrt{n} t^n W_t$, that is, $\widehat{g}_n(s,t)= \sqrt{n} t^n \mathbf{1}_{[0,t]}(s)$. In this case, 
\[
 (\widehat{g}_n\otimes _{1}\widehat{g}_n)(s,t)= n s^nt^n (s\wedge t),
 \]
 which converges to zero weakly in $L^2(\Omega)$, and
 \[
  \langle u_n, DF_{n} \rangle_{\EuFrak H}
  =\int_0^1n t^{2n} W_t^2 dt + n \int_0^1 t^n W_t \left( \int_0^t s^n W_s dW_s \right) dt,
  \]
  which converges in $L^2(\Omega)$ to  $\frac 12 W_1^2$. Therefore, conditions (a) and (b) of 
Proposition \ref{prop-2} are satisfied with $S^2=\frac 12  W_1^2$, and $F_n$ converges in distribution to 
$\frac 1{\sqrt{2}} W_1\times N$, with $N\sim N(0,1)$.  
   One easily see on this particular example that 
the choice  $\widehat{g}_n  =g_n $ does not allows us to conclude in general
(except when $S^2$ is deterministic); indeed,  
one can check here that $\langle u_n, DF_n \rangle_\EuFrak H=
 \frac 1m \| DF_n \|^2 _\EuFrak H$ does not converge in $L^1(\Omega)$.
}
\end{example}

If  we take  $\widehat{g}_n  =g_n $ and  $S^{2}=1$, then condition (b) coincides with condition (iii) in the
introduction. In this case, Nualart and Peccati criterion combined with Lemma
6 in \cite{NO} tells us that, if the sequence of variances converges to one, then condition (a)
is automatically satisfied.

\medskip On the other hand, we can also apply Theorem \ref{main-thm} with $%
u_{n}=g_{n} $. In this way, applying Corollary \ref{cor-1}, we obtain that
the following conditions imply that $F_{n}$ converges to a normal random
variable $N(0,1)$ independent of $X$:

\begin{enumerate}
\item[($\protect\alpha$)] $g_{n}$ converges weakly to zero;

\item[($\protect\beta$)] $\| g_{n}\otimes _{l }g_{n}\|_{\EuFrak %
H^{\otimes 2(q-l )}}$ converges to zero for all $l =1,\dots ,q-1$;

\item[($\protect\gamma$)] $q!\Vert g_{n}\Vert _{\EuFrak H^{\otimes q}}^{2}$
converges to 1.
\end{enumerate}

Indeed, notice first that if $g_{n}$ is bounded in $\EuFrak H^{\odot q}$, then $%
F_{n}$ is bounded in all the Sobolev spaces $\mathbb{D}^{q,p}$, $p\geq 2$.
Then, condition (ii) in Corollary \ref{cor-1} follows from ($\gamma $) and
the equality $D^{q}\left( I_{q}(g_{n})\right) =q!g_{n}$. On the other hand,
condition (i') in Corollary \ref{cor-1} follows from (ii) and%
\begin{eqnarray*}
E\left[ \left\| g_{n}\otimes _{l }D^{l }F_{n}\right\| _{\EuFrak %
H^{\otimes (q-l )}}^{2}\right] &=&\frac{q!^{2}}{(q-l )!^{2}}E\left[
\left\| g_{n}\otimes _{l }I_{q-l }(g_{n})\right\| _{\EuFrak H^{\otimes
(q-l )}}^{2}\right] \\
&=&\frac{q!^{2}}{(q-l )!^{2}}E\left[ \left\| I_{q-l }(g_{n}\otimes
_{l }g_{n})\right\| _{\EuFrak H^{\otimes (q-l )}}^{2}\right] \\
&=&\frac{q!^{2}}{(q-l )!}\left\| g_{n}\widetilde{\otimes}_{l }g_{n}\right\| _{%
\EuFrak H^{\otimes 2(q-l )}}^{2}\\
&\leq &\frac{q!^{2}}{(q-l )!}\left\| g_{n}\otimes_{l }g_{n}\right\| _{%
\EuFrak H^{\otimes 2(q-l )}}^{2}.
\end{eqnarray*}%
In this way we recover the fact that condition (iii) in the introduction
implies the normal convergence.

\setcounter{equation}{0}

\section{Weighted Hermite variations of the fractional Brownian motion}

\label{fbm}

\subsection{Description of the results}

The fractional Brownian motion (fBm) with Hurst parameter $H\in (0,1)$
is a centered Gaussian process $B=\{B_{t},\,t\geq 0\}$ with the covariance
function
\begin{equation}
E(B_{s}B_{t})=R_{H}(s,t)=\frac{1}{2}\left( t^{2H}+s^{2H}-|t-s|^{2H}\right) .
\label{eq1}
\end{equation}%
From (\ref{eq1}), it follows that $E|B_{t}-B_{s}|^{2}=(t-s)^{2H}$ for all $%
0\leqslant s<t$ and that, for each $a>0$, the process $%
\{a^{-H}B_{at},\,t\geq 0\}$ is also a fBm with Hurst parameter $H$ 
(self-similarity property).
As a consequence, the sequence $\{B_{j}-B_{j-1},\,j=1,2,\ldots \}$ is
stationary, Gaussian and ergodic, with correlation given by%
\begin{equation}
\rho _{H}(n)=\frac{1}{2}\left[ |n+1|^{2H}-2|n|^{2H}+|n-1|^{2H}\right],
\label{eq2}
\end{equation}%
which behaves as $H(2H-1)|n|^{2H-2}$ as $n$ tends to infnity.

Set $\Delta B_{k/n}=B_{(k+1)/n}-B_{k/n}$, where $k=0,1,\ldots ,n$, and $n\geq
1 $. The ergodic theorem combined with the self-similarity property implies that
the sequence $n^{2H-1}\sum_{k=0}^{n-1}\left( \Delta
B_{k/n}\right)^{2} $ converges, almost surely and in $L^{1}(\Omega)$, to 
$E(B_{1}^{2})=1$. Moreover, it is well-known (see, e.g., \cite{BM}) that, provided $H\in(0,\frac34)$, 
a central limit theorem holds: the sequence%
\begin{equation}
\frac{1}{\sqrt{n}}\sum_{k=0}^{n-1}\left( n^{2H}\left( \Delta B_{k/n}\right)
^{2}-1\right)
=
\frac{1}{\sqrt{n}}\sum_{k=0}^{n-1}H_2\left(n^{H} \Delta B_{k/n}\right)
  \label{eq2a}
\end{equation}%
converges in law to $N(0,\sigma _{H}^{2})$ as $n\to\infty$, for some constant 
$\sigma_{H}>0$. (Notice also that, by normalizing with 
$\sqrt{n\log n}$ instead of $\sqrt{n}$,
the central limit theorem continues to hold in the critical case $H=\frac34$.)
When $H>\frac34$, the situation is very different. Indeed,  we have in contrast that  
$$
n^{1-2H}\sum_{k=0}^{n-1}\left( n^{2H}\left( \Delta B_{k/n}\right)
^{2}-1\right)
=
n^{1-2H}\sum_{k=0}^{n-1}H_2\left(n^{H} \Delta B_{k/n}\right)
$$
converges in $L^2(\Omega)$.
More generally, consider an integer $q\geq 2$. If $H<1-\frac1{2q}$, then the sequence
\begin{equation}
\frac{1}{\sqrt{n}}\sum_{k=0}^{n-1}H_q\left(n^{H} \Delta B_{k/n}\right)
  \label{eq17a}
\end{equation}%
converges in law to $N(0,\sigma _{q,H}^{2})$ (for some constant $\sigma_{q,H}>0$), whereas,
if $H>1-\frac1{2q}$, then the sequence 
$$
n^{q-qH-1}\sum_{k=0}^{n-1}H_q\left(n^{H} \Delta B_{k/n}\right)
$$
converges in $L^2(\Omega)$.  

Some unexpected results happen 
when we introduce a weight of the form 
$f(B_{k/n})$ in (\ref{eq17a}).
In fact, a new critical value ($H=\frac{1}{2q}$) plays an
important role. 
More precisely, consider the following sequence
of random variables:
\begin{equation}
G_{n}=\frac{1}{\sqrt{n}}\sum_{k=0}^{n-1}f(B_{k/n})H_{q} ( n^{ H} 
\Delta B_{k/n}   ).\label{Gn}
\end{equation}%
Here, the integer $q\geq 2$ is fixed and
the function $f:\mathbb{R}\to\mathbb{R}$ is supposed to satisfy some suitable regularity 
and growth conditions.
In \cite{No1,NNT}, the following convergences as $n\to\infty$ are shown:
\begin{itemize}
\item If $H<\frac1{2q}$, then
\begin{equation}
\ n^{qH-\frac{1}{2}}\,G_{n}\overset{L^{2}(\Omega)}{\longrightarrow }\frac{\left(
-1\right) ^{q}}{2^{q}q!}\int_{0}^{1}f^{(q)}(B_{s})ds.  \label{eq4}
\end{equation}%
\item If $\frac1{2q}<H<1-\frac1{2q}$, then
\begin{equation}
G_{n}\overset{\mathrm{stably}}{\longrightarrow }\,\sigma
_{H,q}\int_{0}^{1}f(B_{s})dW_{s},   \label{eq3}
\end{equation}%
where $W$ is a Brownian motion independent of $B$, and%
\begin{equation}
\sigma _{H,q}^{2}=q!\sum_{r\in \mathbb{Z}}\rho _{H}(r)^{q}<\infty.
\label{eq3a}
\end{equation}
\item If $H=1-\frac1{2q}$, then
$$
\frac{G_n}{\sqrt{\log n}}\,\overset{\mathrm{stably}}{\longrightarrow }\,
\sqrt{\frac2{q!}}
\left(1-\frac1{2q}\right)^{q/2}\left(1-\frac1q\right)^{q/2}
\int_{0}^{1}f(B_{s})dW_{s},   
$$
where $W$ is a Brownian motion independent of $B$.
\item If $H>1-\frac1{2q}$, then
$$
n^{q(1-H)-\frac12}\,G_n\,\overset{L^2(\Omega)}{\longrightarrow }\,
\int_{0}^{1}f(B_{s})dZ^{(q)}_{s},   
$$
where $Z^{(q)}$ denotes the Hermite process of order $q$ canonically constructed from $B$
(see \cite{NNT} for the details).
\end{itemize}
In addition, when $q=2$ and $H=\frac14$, it was shown in \cite{NR} that $G_n$ converges stably 
to a linear combination of the limits in (\ref%
{eq3}) and (\ref{eq4}). (The proof of this last result follows an approach similar
to the proof of our Theorem \ref{main-thm}, and allows to derive  
a change of variable formula for the fBm of Hurst index $\frac14$, 
with a correction term that is an ordinary It\^o integral with respect to a Brownian motion 
that is independent of $B$.)
But the convergence of $G_n$ in the critical case $H=\frac1{2q}$, $q\geq 3$, was open till now.

\medskip In the present paper, we are going to show that Theorem \ref{main-thm} provides a proof
of the following new result, valid for any integer $q\geq 2$ and any index 
$H\in \left( \frac{1}{4q},\frac{1}{2}\right)$: 
\begin{equation}
G_{n}-n^{-\frac{1}{2}-qH}\frac{\left( -1\right) ^{q}}{2^{q}q!}
\sum_{k=0}^{n-1} f^{(q)}(B_{k/n})
\overset{\mathrm{stably}}{\longrightarrow }%
\,\sigma _{H,q}\int_{0}^{1}f(B_{s})dW_{s}.  \label{eq11}
\end{equation}
(See Theorem \ref{thm-eq11} below for a precise statement.)
Notice that (\ref{eq11}) provides a new proof of (\ref{eq3})
in the case $H\in \left( \frac{1}{2q},\frac{1}{2}\right)$
(without considering two different levels of discretization $%
n\leqslant m$, as in \cite{NNT}).
More importantly, in the critical case $H=\frac1{2q}$, convergence (\ref{eq11}) yields:
$$
G_{n}
\overset{\mathrm{stably}}{\longrightarrow }%
\,
\frac{\left( -1\right) ^{q}}{2^{q}q!}%
\int_{0}^{1}f^{(q)}(B_{s})ds+
\sigma _{1/(2q),q}\int_{0}^{1}f(B_{s})dW_{s}. 
$$
Hence, the understanding of the asymptotic behavior of the weighted Hermite variations of the fBm is now complete
(indeed, the case $H=\frac{1}{2q}$, $q\geq 3$, was the only remaining case, as mentioned in the discussion
above).

The main idea of the proof of (\ref{eq11}) is a decomposition of the random variable $G_n$
  using equation (\ref{t3}). The term with $r=0$
 is a multiple Skorohod integral of order $q$ 
and, by Theorem \ref{thm1} below, it converges in law for any $H\in \left( \frac{1}{4q},%
\frac{1}{2}\right)$.
The  term with $r=q$  behaves as $%
-n^{-\frac12-qH}\frac{\left( -1\right) ^{q}}{2^{q}q!}\sum_{k=0}^{n}f^{(q)}(B_{k/n})$. 
The remaining terms ($1\le r \le q-1$) converge to zero in $L^2(\Omega)$.

\subsection{Some preliminaries on the fractional Brownian motion}

Before proving (\ref{eq11}), we need some preliminaries on the
Malliavin calculus associated with the fBm and some technical results (see %
\cite[Chapter 5]{Nbook}).

In the following we assume $H\in \left( 0,\frac{1}{2}\right) $. We denote by
$\mathcal{E}$ the set of step functions on $[0,1]$. Let $\EuFrak H$ be the
Hilbert space defined as the closure of $\mathcal{E}$ with respect to the
scalar product
\begin{equation*}
\left\langle {\mathbf{1}}_{[0,t]},{\mathbf{1}}_{[0,s]}\right\rangle _{%
\EuFrak H}=R_{H}(t,s)=\frac12\big(s^H+t^H-|t-s|^H\big).
\end{equation*}%
The mapping $\mathbf{1}_{[0,t]}\rightarrow B_{t}$ can be extended to a
linear isometry between the Hilbert space $\EuFrak H$ and the Gaussian space
spanned by $B$. 
We denote this isometry by $\phi \rightarrow B(\phi )$. In
this way $\{B(\phi ),\,\phi \in \EuFrak H\}$ is an isonormal Gaussian space.
(In fact, we know that the space $\EuFrak H$ coincides with $I_{0+}^{H-\frac{1}{2}%
}(L^{2}[0,1])$, where $$I_{0+}^{H-\frac{1}{2}}f(x)=\frac{1}{\Gamma (H-\frac{1%
}{2})}\int_{0}^{x}(x-y)^{H-\frac{3}{2}}f(y)dy$$ 
is the left-sided Liouville
fractional integral of order $H-\frac{1}{2}$, see \cite{DU}.)

From now on, we will make use of the notation%
\begin{eqnarray*}
&&\varepsilon _{t}=\mathbf{1}_{[0,t]}, \\
&&\partial _{k/n}= \varepsilon _{(k+1)/n}-\varepsilon _{k/n}=\mathbf{1}%
_{(k/n,(k+1)/n]},  \label{eq5}
\end{eqnarray*}%
for $t\in \lbrack 0,1]$, $n\geq 1$, and $k=0,\ldots ,n-1$. Notice that $%
H_q(n^H\Delta B_{k/n}) =n^{qH}I_{q}(\partial _{k/n}^{\otimes q})$.

\medskip We need the following technical lemma.

\begin{lemma}
\label{lem1} Recall that $H<\frac{1}{2}$. Let $n\geq 1$ and $k=0,\ldots
,n-1 $. We have

\begin{itemize}
\item[(a)] $\big|E(B_{r}(B_{t}-B_{s}))\big|\leqslant (t-s)^{2H}$ for any $\ r\in $ $%
[0,1]$ and $0\leqslant s<t\leqslant 1$.

\item[(b)] $\left| \left\langle \varepsilon _{t},\;\partial
_{k/n}\right\rangle _{\EuFrak H}\right| \leqslant n^{-2H}$ for any $t\in
\lbrack 0,1]$.

\item[(c)] $\sup_{t\in \lbrack 0,1]}\sum_{k=0}^{n-1}\left| \left\langle
\varepsilon _{t},\;\partial _{k/n}\right\rangle _{\EuFrak H}\right| =O(1)$
as $n$ tends to infinity.

\item[(d)] For any integer $q\geq 2$,
\begin{equation}
\sum_{k=0}^{n-1}\left| \left\langle \varepsilon _{k/n},\;\partial
_{k/n}\right\rangle _{\EuFrak H}^{q}-\frac{(-1)^{q}}{2^{q}n^{2qH}}\right|
=O(n^{-2H(q-1)})\quad\mbox{as $n$ tends to infinity}.  \label{d1}
\end{equation}

\item[(e)] Recall the definition (\ref{eq2}) of $\rho_H$. We have $$\left\langle \partial _{j/n},\;\partial
_{k/n}\right\rangle _{\EuFrak H}=n^{-2H}\,\rho _{H}(k-j).$$ Consequently, for
any integer $q\geq 1$, we can write%
\begin{equation}
\sum_{k,j=0}^{n-1}\left| \left\langle \partial _{j/n},\;\partial
_{k/n}\right\rangle _{\EuFrak H}\right| ^{q}=O(n^{1-2qH})\quad\mbox{as $n$ tends to infinity}.  \label{d3a}
\end{equation}
\end{itemize}
\end{lemma}

\begin{proof}
We have%
\begin{eqnarray*}
E(B_{r}(B_{t}-B_{s})) &=&\frac{1}{2}\left( r^{2H}+t^{2H}-|t-r|^{2H}\right) -%
\frac{1}{2}\left( r^{2H}+s^{2H}-|s-r|^{2H}\right) \\
&=&\frac{1}{2}(t^{2H}-s^{2H})+\frac{1}{2}\left( |s-r|^{2H}-|t-r|^{2H}\right).
\end{eqnarray*}%
Using the inequality $|b^{2H}-a^{2H}|\leqslant |b-a|^{2H}$ for any $a,b\in
\lbrack 0,1]$, we deduce (a). Property (b) is an
immediate consequence of (a). To show property (c) we use%
\begin{equation*}
\left\langle \varepsilon _{t},\;\partial _{k/n}\right\rangle _{\EuFrak H}=%
\frac{1}{2n^{2H}}\left[ (k+1)^{2H}-k^{2H}-|k+1-nt|^{2H}+|k-nt|^{2H}\right] .
\end{equation*}%
Property (d) follows from%
\begin{equation*}
\left\langle \varepsilon _{k/n},\;\partial _{k/n}\right\rangle _{\EuFrak H}=%
\frac{1}{2n^{2H}}\left[ (k+1)^{2H}-k^{2H}-1\right] ,
\end{equation*}%
and%
\begin{eqnarray*}
\left| \left\langle \varepsilon _{k/n},\;\partial _{k/n}\right\rangle _{%
\EuFrak H}^{q}-\frac{(-1)^{q}}{2^{q}n^{2qH}}\right|  &=&\frac{1}{2^{q}n^{2qH}%
}\left| \left[ (k+1)^{2H}-k^{2H}-1\right] ^{q}-(-1)^{q}\right|  \\
&= &\frac{1}{2^{q}n^{2qH}}\sum_{i=1}^{q}\binom{q}{i}\left[
(k+1)^{2H}-k^{2H}\right] ^{i}\\
&\leq &\frac{1}{2^{q}n^{2qH}}\left[
(k+1)^{2H}-k^{2H}\right]\sum_{i=1}^{q}\binom{q}{i}.
\end{eqnarray*}

Finally, property (e) follows from%
\begin{equation*}
\sum_{k,j=0}^{n-1}\left| \left\langle \partial _{j/n},\;\partial
_{k/n}\right\rangle _{\EuFrak H}\right|  ^q\leq n^{-2qH}\sum_{k,j=0}^{n-1}\left|
\rho _{H}(j-k)\right| ^q
\leq n^{1-2qH}\sum_{r\in\mathbb{Z}}|\rho_H(r)|^q.
\end{equation*}
\end{proof}

\subsection{An auxiliary convergence result}
 
From now on, we fix $q\geq 2$ and we make use of the following hypothesis on
$f:\mathbb{R}\rightarrow \mathbb{R}$:

\medskip \textbf{(H)} $f$ belongs to $\mathcal{C}^{2q}$ and,
for any $p\geq 2$ and $i=0,\dots, 2q$,
\begin{equation}
E(\sup_{t\in \lbrack 0,1]}|f^{(i)}(B_{t})|^{p})<\infty .  \label{n11}
\end{equation}%

Notice that a sufficient condition for (\ref{n11}) to hold is that $f$
satisfies an exponential growth condition of the form $\left|
f^{(2q)}(x)\right| \leqslant ke^{c|x|^{p}}$ for some constants $c,k>0$ and $%
0<p<2$.

\bigskip The aim of this section is to prove the following auxiliary convergence result.

\begin{theorem}
\label{thm1}Suppose $H\in \left( \frac{1}{4q},\frac{1}{2}\right) $, and let $%
f$ be a function satisfying Hypothesis \textbf{(H)}. Consider the sequence
of $q$-parameter step processes defined by%
\begin{equation}
u_{n}=n^{qH-\frac{1}{2}} \sum_{k=0}^{n-1}f(B_{k/n})\partial _{k/n}^{\otimes
q}.  \label{v4}
\end{equation}%
Then $u_n\in{\rm Dom}\delta^q$, and $\delta ^{q}(u_{n})$ converges stably to $\sigma
_{H,q}\int_{0}^{1}f(B_{s})dW_{s}$,   where $W$ is a Brownian motion
independent of $B$, and  $\sigma_{H,q}>0$ is defined in (\ref{eq3a}).
\end{theorem}

\begin{proof}
The fact that $u_n$ belongs to ${\rm Dom}\delta^q$ is a consequence of the inclusion
$\mathbb{D}^{q,2}(\EuFrak H^{\otimes q})\subset {\rm Dom}\delta^q$ and hypothesis 
\textbf{(H)}.
We are now going to show that the sequence $F_{n}= \delta ^{q}(u_{n})$
satisfies the conditions of Theorem \ref{main-thm}.  We make use
of the notation
\begin{equation}\label{star}
\alpha _{k,j}=\left\langle \varepsilon _{k/n},\partial _{j/n}\right\rangle _{%
\EuFrak H},\;\beta _{k,j}=\left\langle \partial _{k/n},\partial
_{j/n}\right\rangle _{\EuFrak H},
\end{equation}%
for $k,j=0,\ldots ,n-1$ and $n\geq 1$. Also $C$ will denote a generic
constant.

\medskip \textit{Step 1. } Let us show first that $F_{n}$ is bounded in $%
L^{2}(\Omega)$. Taking into account the continuity of the Skorohod integral from
the space $\mathbb{D}^{q,2}(\EuFrak H^{\otimes q})$ into $L^{2}(\Omega)$ (see (\ref%
{Me2})), it suffices to show that $\ u_{n}$ is bounded in $\mathbb{D}^{q,2}(%
\EuFrak H^{\otimes q})$. Actually we are going to show that $u_{n}$ is
bounded in $\mathbb{D}^{k,p}(\EuFrak H^{\otimes k})$ for any integer $%
k\leq  2q$ and any real number $p\geq 2$.  Using the estimate (\ref{d3a}) we obtain 
\begin{equation*}
\left\| u_{n}\right\| _{\EuFrak H^{\otimes q}}^{2}=n^{2qH-1}\
\sum_{k,j=0}^{n-1}f(B_{k/n})f(B_{j/n})\beta _{k,j}^{q}\leq C\sup_{0\leq
t\leq 1}\left| f(B_{t})\right| ^{2}.
\end{equation*}%
Moreover for any integer $k\geq 1$, 
\begin{equation*}
D^{k}u_{n}=n^{qH-\frac{1}{2}}\sum_{j=0}^{n-1}f^{(k)}(B_{j/n})\varepsilon
_{j/n}^{\otimes k}\otimes \partial _{j/n}^{\otimes q},
\end{equation*}%
and we obtain in the same way%
\begin{eqnarray*}
\left\| D^{k}u_{n}\right\| _{\EuFrak H^{\otimes (q+k)}}^{2} &=&n^{2qH-1}\
\sum_{l,j=0}^{n-1}f^{(k)}(B_{l/n})f^{(k)}(B_{j/n})\left\langle \varepsilon
_{l/n},\varepsilon _{j/n}\right\rangle ^{k}\beta _{l,j}^{q} \\
&\leq &C\sup_{0\leq t\leq 1}\left| f^{(k)}(B_{t})\right| ^{2}.
\end{eqnarray*}%
Then the result follows from hypothesis \textbf{(H)}.

 \medskip \textit{Step 2. }Let us show condition (i) of Theorem \ref{main-thm}. 
Fix some integers $r,k_{1},\ldots ,k_{q-1}\geq 0$ such that $%
k_{1}+2k_{2}+\ldots +(q-1)k_{q-1}+r=q$. Let $h\in \EuFrak H^{\otimes r}$. We
claim that $\left\langle u_{n},(DF_{n})^{\otimes k_{1}}\otimes \ldots
\otimes (D^{q-1}F_{n})^{\otimes k_{q-1}}\otimes h\right\rangle _{\EuFrak %
H^{\otimes q}}$ converges to zero in $L^{1}(\Omega)$. Suppose first that $r\geq 1$.
Without loss of generality, we can assume that $h$ has the form $g \otimes \varepsilon _{t}$,
with $g\in\EuFrak H^{\otimes(r-1)}$. Set $%
\Phi _{n}=(DF_{n})^{\otimes k_{1}}\otimes \ldots \otimes
(D^{q-1}F_{n})^{\otimes k_{q-1}}\otimes g$. Then we can write%
\begin{equation*}
\left\langle u_{n}, \Phi _{n}\otimes \varepsilon_t\right\rangle _{\EuFrak H^{\otimes
q}}=n^{qH-\frac{1}{2}} \sum_{k=0}^{n-1}f(B_{k/n})\left\langle \partial
_{k/n}^{\otimes (q-1)},\Phi _{n}\right\rangle _{\EuFrak H^{\otimes
(q-1)}}\left\langle \partial _{k/n},\varepsilon _{t}\right\rangle _{\EuFrak %
H }.
\end{equation*}%
As a consequence,%
\begin{eqnarray*}
E\left( \left| \left\langle u_{n},  \Phi _{n}\otimes \varepsilon_t\right\rangle _{\EuFrak
H^{\otimes q}}\right| \right) &\leq &n^{qH-\frac{1}{2}}\
\sum_{k=0}^{n-1}E\left( \left| f(B_{k/n})\left\langle \partial
_{k/n}^{\otimes (q-1)},\Phi _{n}\right\rangle _{\EuFrak H^{\otimes
(q-1)}}\right| \right)  \\
&& \times \left| \left\langle \partial _{k/n},\varepsilon
_{t}\right\rangle _{\EuFrak H}\right| .
\end{eqnarray*}%
Condition (c) of Lemma \ref{lem1} implies%
\begin{equation*}
\sum_{k=0}^{n-1}\left| \left\langle \partial _{k/n},\varepsilon
_{t}\right\rangle _{\EuFrak H}\right| \leq C.
\end{equation*}%
Hence,%
\begin{equation*}
E\left( \left| \left\langle u_{n}, \Phi _{n}\otimes \varepsilon_t\right\rangle _{\EuFrak
H^{\otimes q}}\right| \right) \leq Cn^{H-\frac{1}{2}}\left( E\left( \left\|
\Phi _{n}\right\| _{\EuFrak H^{\otimes (q-1)}}^{2}\right) \right) ^{\frac{1}{%
2}}.
\end{equation*}%
On the other hand%
\begin{equation*}
\left\| \Phi _{n}\right\| _{\EuFrak H^{\otimes (q-1)}}^{2}=\left\| g\right\|
_{\EuFrak H^{\otimes (r-1)}}^{2}\prod_{m=1}^{q-1}\left\| D^{m}F_{n}\right\|
_{\EuFrak H^{\otimes m}}^{2k_{m}},
\end{equation*}%
and applying the generalized H\"{o}lder's inequality%
\begin{eqnarray*}
E\left( \left\| \Phi _{n}\right\| _{\EuFrak H^{\otimes (q-1)}}^{2}\right) 
&\leq &C\prod_{m=1}^{q-1}\left( E\left( \left\| D^{m}F_{n}\right\| _{\EuFrak %
H^{\otimes m}}^{2k_{m}(q-1)}\right) \right) ^{\frac{1}{q-1}} \\
&=&C\prod_{m=1}^{q-1} \left\| D^{m}F_{n}\right\| _{L^{2k_{m}(q-1)}(\Omega ;%
\EuFrak H^{\otimes m})}^{2k_{m} }.
\end{eqnarray*}%
By Meyer's inequalities (\ref{Me2}), for any $1\leq m\leq q-1$ and any $p\geq 2$, we obtain, using Step 1, that%
\begin{eqnarray*}
\left\| D^{m}F_{n}\right\| _{L^{p}(\Omega ;\EuFrak H^{\otimes m})}^{\
}&=&\left\| D^{m}\delta ^{q}(u_{n})\right\| _{L^{p}(\Omega ;\EuFrak H^{\otimes
m})}\\
&\leq& 
\left\| \delta^q(u_n)\right\|_{\mathbb{D}^{m,p}} \leq
C\left\| u_{n}\right\| _{\mathbb{D}^{m+q,p}(\EuFrak H^{\otimes
q})}\leq C.
\end{eqnarray*}
Therefore,%
\begin{equation*}
E\left( \left| \left\langle u_{n}, \Phi _{n}\otimes \varepsilon_t\right\rangle _{\EuFrak
H^{\otimes q}}\right| \right) \leq Cn^{H-\frac{1}{2}},
\end{equation*}%
which converges to zero as $n$ tends to infinity because $H<\frac{1}{2}$.

Suppose now that $r=0$. In this case, we have $\Phi _{n}=(DF_{n})^{\otimes k_{1}}\otimes
\cdots \otimes (D^{q-1}F_{n})^{\otimes k_{q-1}}$. \ Then%
\begin{equation}
\left\langle \partial _{j/n}^{\otimes q},\Phi _{n}\right\rangle _{\EuFrak %
H^{\otimes q}}=\left\langle \partial _{j/n},DF_{n}\right\rangle _{\EuFrak %
H}^{k_{1}}\cdots \left\langle \partial _{j/n}^{\otimes
(q-1)},D^{q-1}F_{n}\right\rangle _{\EuFrak H^{\otimes (q-1)}}^{k_{q-1}}.
\label{t9}
\end{equation}%
From (\ref{t9}) and (\ref{v4}) we obtain%
\begin{equation}  \label{x1}
\left\langle u_{n}, \Phi _{n}\right\rangle _{\EuFrak H^{\otimes q}}=n^{qH-%
\frac{1}{2}} \sum_{k=0}^{n-1}f(B_{k/n})\prod_{m=1}^{q-1}\left\langle
\partial _{j/n}^{\otimes m},D^{m}F_{n}\right\rangle _{\EuFrak H^{\otimes m}}^{k_{m}}.
\end{equation}%
Notice that \ for any $m=1,\ldots ,q-1$, the term $\left\langle \partial
_{j/n}^{\otimes m},D^{m}F_{n}\right\rangle _{\EuFrak H^{\otimes m}}$ can be
estimated by $n^{-mH}\left\| D^{m}F_{n}\right\| _{\EuFrak H^{\otimes m}}$.
Then, taking into account that $$\sup_{n}E\left( \left\| D^{m}F_{n}\right\| _{%
\EuFrak H^{\otimes m}}^{p}\right) <\infty $$ for any $p\geq 2$, and that $%
\sum_{m=1}^{q-1}mk_{m}=q$, we obtain for  $E\left( \left|
\left\langle u_{n}, \Phi _{n}\right\rangle _{\EuFrak H^{\otimes q}}\right|
\right) $ an estimate of the form $C\sqrt{n}$, which is unfortunately not satisfactory. For
this reason, a finer analysis of the terms $\left\langle \partial
_{j/n}^{\otimes m},D^{m}F_{n}\right\rangle _{\EuFrak H^{\otimes m}}$ is
required.

First we are going to apply formula (\ref{t5}) to compute the derivative $%
D^{m}F_{n}$, $m=1,\ldots ,q-1$:
\begin{eqnarray}  \notag
D^{m}F_{n} &=&\sum_{i=0}^{m}\binom{m}{i}\binom{q}{i}i!\delta
^{q-i}(D^{m-i}u_{n}) \\  \notag
&=&n^{qH-\frac{1}{2}}\sum_{i=0}^{m}\binom{m}{i}\binom{q}{i}i!\
\sum_{l=0}^{n-1}\left( \varepsilon _{l/n}^{\otimes (m-i)}\otimes \partial
_{l/n}^{\otimes i}\right)  \\  \label{x2}
&&\times \delta ^{q-i}\left( f^{(m-i)}(B_{l/n})\partial
_{l/n}^{\otimes (q-i)}\right).
\end{eqnarray}%
Set  $ \Psi  _{n}^{m,j} = \left\langle \partial _{j/n}^{\otimes m},D^{m}F_{n}\right\rangle
_{\EuFrak H^{\otimes m}}$, and recall the definition of $\alpha_{k,j}$
and $\beta_{k,j}$ from (\ref{star}).  From  (\ref{x2}) we obtain
\begin{eqnarray}
\Psi  _{n}^{m,j}&=&n^{qH-\frac{1}{2}}\sum_{i=0}^{m}\binom{m}{i}\binom{q}{i}i!
\sum_{l=0}^{n-1}\alpha _{l,j}^{m-i}\beta _{l,j}^{i}\,\delta ^{q-i}\left(
f^{(m-i)}(B_{l/n})\partial _{l/n}^{\otimes (q-i)}\right)   \notag \\
&=&\sum_{i=0}^{m}\Phi _{n}^{i,m,j},  \label{t10}
\end{eqnarray}%
with
\begin{equation*}
\Phi _{n}^{i,m,j}=n^{qH-\frac{1}{2}}\binom{m}{i}\binom{q}{i}i! \sum_{l=0}^{n-1}\alpha
_{l,j}^{m-i}\beta _{l,j}^{i}\delta ^{q-i}\left( f^{(m-i)}(B_{l/n})\partial
_{l/n}^{\otimes (q-i)}\right) .
\end{equation*}%
By Meyer inequalities (\ref{Me2}) we obtain, using also assumption {\bf (H)}, that, for any $p\geq 2$,%
\begin{eqnarray}  \notag
\left\| \delta ^{q-i}\left( f^{(m-i)}(B_{l/n})\partial _{l/n}^{\otimes
(q-i)}\right) \right\| _{L^p} &\leq &C\left\| f^{(m-i)}(B_{l/n})\partial
_{l/n}^{\otimes (q-i)}\right\| _{\mathbb{D}^{q-i,p}(\EuFrak H^{\otimes q-i})}
\\   \label{x3}
&\leq &Cn^{-(q-i)H}.
\end{eqnarray}%
Using Lemma \ref{lem1} (b) and (e) we have $\left| \alpha
_{l,j}^{m-i}\right| \leq Cn^{-(m-i)2H}$ and $\sum_{l=0}^{n-1}\left| \beta
_{l,j}^{i}\right| \leq Cn^{-2iH}$. Therefore, for any $i\geq 1$, we have
\begin{equation}
\left\| \Phi _{n}^{i,m,j}\right\| _{L^p}\leq Cn^{iH-\frac{1}{2}}\
\sum_{l=0}^{n-1}\left| \alpha _{l,j}^{m-i}\beta _{l,j}^{i}\right| \leq Cn^{-%
\frac{1}{2}-2mH+iH}.  \label{t12}
\end{equation} 
On the other hand, if $i=0$, Lemma    \ref{lem1} (c)  and  (\ref{x3}) yield
\begin{equation}
\left\| \Phi _{n}^{0,m,j}\right\| _{L^p}\leq C   n^{-%
\frac{1}{2}-2mH+2H}.  \label{t12a}
\end{equation} 
Notice that the estimate for the $L^p(\Omega)$-norm 
of  $\Phi _{n}^{0,m,j}$ in the case $i=0$ is worst than for $i\ge 1$. We will see later that, for $p=2$, 
we can get a better estimate for $\Phi _{n}^{0,m,j}$.

Because $\sum_{m=1}^{q-1}k_{m}\geq 2$,  the number of factors
in $\prod_{m=1}^{q-1}\left\langle \partial
_{j/n},D^{m}F_{n}\right\rangle _{\EuFrak H^{\otimes m}}^{k_{m}}$ is at least
two. As a consequence, we can write%
\begin{equation*}
\left\langle u_{n}, \Phi _{n}\right\rangle _{\EuFrak H^{\otimes q}}=n^{qH-%
\frac{1}{2}}\sum_{j=0}^{n-1}f(B_{j/n})\Psi _{n}^{\mu ,j}\Psi _{n}^{\nu
,j}\Theta _{n}^{j},
\end{equation*}%
for some $\mu $, $\nu $ (not necessarily distinct), where 
\begin{equation}  \label{x4}
\Theta _{n}^{j}=\left( \Psi _{n}^{\mu ,j}\right) ^{k_{\mu
}-1}\left( \Psi _{n}^{\nu ,j}\right) ^{k_{\nu }-1}\prod_{\substack{ m=1 \\ m\neq \mu ,\nu }}^{q-1} \left(
\Psi _{n}^{m,j}\right) ^{k_{m}}.
\end{equation}%
Consider the decomposition%
\begin{equation*}
\left\langle u_{n}, \Phi _{n}\right\rangle _{\EuFrak H^{\otimes
q}}=A_{n}+ B_n,
\end{equation*}%
where%
\begin{eqnarray*}
A_{n} &=&n^{qH-\frac{1}{2}}\sum_{j=0}^{n-1}f(B_{j/n})\left(
\sum_{i=0}^{\mu} \sum_{k=0}^{\nu}  \mathbf{1}_{i+k \ge 1} \Phi _{n}^{i,\mu ,j} \Phi
_{n}^{k,\nu ,j}\right)  \Theta _{n}^{j}, \\
B_{n} &=& n^{qH-\frac{1}{2}}\sum_{j=0}^{n-1}f(B_{j/n})\Phi _{n}^{0,\mu
,j}\Phi _{n}^{0,\nu ,j} \Theta _{n}^{j}.
\end{eqnarray*}%
From (\ref{x4}) and the estimate     $ \| \Psi^{m,j}_n \| _{L^p} \le Cn^{-mH}$,  for all $p\geq 2$ and $1\le m\le q$, we obtain
\begin{equation}
\left\| \Theta _{n}^{j}\right\| _{L^p}\leq Cn^{-H(q-\mu -\nu )}.  \label{t11}
\end{equation}%
Then, from (\ref{t12}), (\ref{t12a}) and (\ref{t11}) we obtain%
\begin{eqnarray*}
E\left( \left| A_{n}\right| \right)  &\leq &Cn^{qH+\frac{1}{2}}n^{-H(q-\mu
-\nu )}  \Big( \sum _{i=1}^\mu \sum _{k=1}^\nu  n^{-1-2(\mu +\nu )H+(i+k)H}    \\
&&\quad +\sum _{i=1}^\mu n^{-1-2(\mu +\nu )H+ iH+2H} 
+\sum _{k=1}^\nu n^{-1-2(\mu +\nu )H+ kH+2H}  \Big) \\
&=&C n^{-\frac 12}+   n^{ -\frac 12 +2H-\mu H }+   n^{ -\frac 12 +2H-\nu H } ,
\end{eqnarray*}%
which converges to zero as $n$ tends to infinity, because $\mu, \nu \ge 1$ and $H<\frac 12$.

For the term $B_n$ using again  the estimates  (\ref{t12a}) and (\ref{t11}) we get
\begin{eqnarray*}
E\left( \left| B_{n}\right| \right)   &\leq & C n^{qH + \frac 12 -H(q-\mu-\nu) -1 -2H(\mu+\nu) +4H}
=C n^{ - \frac 12    - H(\mu+\nu) +4H}  \\
&\leq&  C n^{ - \frac 12     +2H},
\end{eqnarray*}
which converges to zero as $n$ tends to infinity if $H<\frac 14$.  
To handle the case $H\in \left[ \frac 14 ,\frac 12 \right)$ we need more precise estimates for the $L^2(\Omega)$-norm of
$\Phi^{0,\nu,j}_n$. 
 We have, using formula (\ref{t13}) 
\begin{eqnarray*}
&&E\left[ \left( \Phi _{n}^{0,\nu ,j}\right) ^{2}\right] =
\binom{q}{i}^2\binom{m}{i}^2i!^2
E\left( \left| n^{qH-\frac{1}{2}} \sum_{l=0}^{n-1}\alpha _{l,j}^{\nu
}\delta ^{q}\left( f^{(\nu )}(B_{l/n})\partial _{l/n}^{\otimes q}\right)
\right| ^{2}\right)  \\
&=&n^{2qH-1}\binom{q}{i}^2\binom{m}{i}^2i!^2\sum_{l,l^{\prime }=0}^{n-1}\alpha _{l,j}^{\nu }\alpha
_{l^{\prime },j}^{\nu }\\
&&\hskip4cm\times
E\left( \delta ^{q}\left( f^{(\nu )}(B_{l/n})\partial
_{l/n}^{\otimes q}\right) \delta ^{q}\left( f^{(\nu )}(B_{l^{\prime
}/n})\partial _{l^{\prime }/n}^{\otimes q}\right) \right)  \\
&=&n^{2qH-1}\binom{q}{i}^2\binom{m}{i}^2i!^2\sum_{l,l^{\prime }=0}^{n-1}\alpha _{l,j}^{\nu }\alpha
_{l^{\prime },j}^{\nu }\sum_{i=0}^{q}\binom{q}{i}^{2}i!\alpha _{l,l^{\prime
}}^{q-i}\alpha _{l^{\prime },l}^{q-i}\beta _{l,l^{\prime }}^{2i}\\
&&\hskip5cm\times E\left(
f^{(\nu +q-i)}(B_{l/n})\ f^{(\nu +q-i)}(B_{l^{\prime }/n}) \right)  \\
&=&\sum_{i=0}^{q}R_{n}^{i}\text{.}
\end{eqnarray*}%
If $i\geq 1$, then  $\sum_{l,l^{\prime }=0}^{n-1}\beta _{l,l^{\prime
}}^{2i}\leq Cn^{1-4iH}$, and we obtain an estimate of the form $\left\|
R_{n}^{i}\right\| _{L^2}\leq Cn^{\gamma }$, where%
\begin{equation*}
\gamma =\frac{1}{2}\left( 2qH-1-4\nu H-4(q-i)H+1-4iH\right) =-qH-2\nu H.
\end{equation*}%
For $i=0$, then $\sup_{n}\sum_{l,l^{\prime }=0}^{n-1}\left| \alpha
_{l,l^{\prime }}\alpha _{l^{\prime },l}\right| <\infty $, and we get%
\begin{equation*}
\gamma =\frac{1}{2}\left( 2qH-1-2H(2\nu +2q-2)\right) =-qH-2\nu H-\frac{1}{2}%
+2H.
\end{equation*}%
We have obtained the estimate%
\begin{equation}
\left\| \Phi _{n}^{0,\nu ,j}\right\| _{L^2}\leq Cn^{-qH-2\nu H+   2H-%
\frac{1}{2} }.  \label{a1}
\end{equation}%
 Fix   $  \frac{1}{4qH}  <\alpha <1$. This choice is possible because $\frac{1}{%
4qH}<1$. We have, by H\"{o}lder's inequality,%
\begin{equation*}
E\left( \left| B_{n}\right| \right) \leq Cn^{qH-\frac{1}{2}%
}\sum_{j=0}^{n-1}\left\| \Phi _{n}^{0,\mu ,j}\right\| _{L^2}^{\alpha }\left\|
  \Phi _{n}^{0,\nu ,j}\right\| _{L^2}^{\alpha }\left\|   \left|   \Phi
_{n}^{0,\mu ,j}\Phi _{n}^{0,\nu ,j}  \right| ^{ 1-\alpha  }\Theta
_{n}^{j}\right\| _{L^{\frac 1 {1-\alpha}}}.
\end{equation*}%
Using   (\ref{a1}),  (\ref{t12a}) and (\ref{t11}) we obtain
\begin{equation}
E\left( \left| B_{n}\right| \right) \leq Cn^{\gamma },  \label{a2}
\end{equation}%
where%
\begin{eqnarray*}
\gamma  &=&qH+\frac{1}{2}+\left[ -2qH-2(\mu +\nu )H+ 4H-1%
 \right]  \alpha  \\
&&-H(q-\mu -\nu )+ (1-\alpha) (-1 -2H(\mu+\nu) +4H) \\
&=& -\frac 12 +4H -H(\mu+\nu) -2\alpha qH  \\
&   \le &  -\frac 12  +2H  -2\alpha qH
\le \frac 12  -2\alpha qH<0,
\end{eqnarray*}
because $H<\frac 12$.
 Therefore $E\left( \left| B_{n}\right| \right) $ converges to zero as $n$
tends to infinity.  

\textit{Step 3. }Let us show condition (ii). \ We have%
\begin{equation*}
\left\langle u_{n},D^{q}F_{n}\right\rangle _{\EuFrak H^{\otimes q}}=n^{qH-%
\frac{1}{2}}\sum_{j=0}^{n-1}f(B_{j/n})\left\langle \partial _{j/n}^{\otimes
q},D^{q}F_{n}\right\rangle _{\EuFrak H^{\otimes q}}.
\end{equation*}%
From (\ref{t10}) we get%
\begin{equation*}
\left\langle \partial _{j/n}^{\otimes q},D^{q}F_{n}\right\rangle _{\EuFrak %
H^{\otimes q}}=n^{qH-\frac{1}{2}}\sum_{i=0}^{q}\binom{q}{i}^2i!\
\sum_{l=0}^{n-1}\alpha _{l,j}^{q-i}\beta _{l,j}^{i}\delta ^{q-i}\left(
f^{(q-i)}(B_{l/n})\partial _{l/n}^{\otimes (q-i)}\right) .
\end{equation*}%
Therefore, we can make the decomposition%
\begin{equation*}
\left\langle u_{n},D^{q}F_{n}\right\rangle _{\EuFrak H^{\otimes
q}}=A_{n}+B_{n}+C_{n},
\end{equation*}%
where 
\begin{eqnarray*}
A_{n} &=&n^{2qH-1}q!\sum_{l,j=0}^{n-1}\beta _{l,j}^{q}f(B_{l/n})f(B_{j/n}),
\\
B_{n} &=&n^{2qH-1}\sum_{i=1}^{q-1}\binom{q}{i}^2i! \sum_{l,j=0}^{n-1}\alpha
_{l,j}^{q-i}\beta _{l,j}^{i}f(B_{j/n})\delta ^{q-i}\left(
f^{(q-i)}(B_{l/n})\partial _{l/n}^{\otimes (q-i)}\right) , \\
C_{n} &=&n^{2qH-1} \sum_{l,j=0}^{n-1}\alpha _{l,j}^{q}f(B_{j/n})\delta
^{q}\left( f^{(q)}(B_{l/n})\partial _{l/n}^{\otimes (q)}\right) .
\end{eqnarray*}%
The term $A_{n}$ converges to a nonnegative square integrable random
variable. Indeed,%
\begin{eqnarray*}
A_{n} &=&\frac{q!}{2^{q}n}\sum_{k,j=0}^{n-1}f(B_{k/n})f(B_{j/n})\left(
|k-j+1|^{2H}+|k-j-1|^{2H}-2|k-j|^{2H}\right) ^{q} \\
&=&\frac{q!}{2^{q}n}\sum_{p=-\infty }^{\infty }
\!\!\!\!\sum_{j=0\vee -p}^{\left(
n-1\right) \wedge \left( n-1-p\right) }  \!\!\!\!\!\!\!\!f(B_{j/n})f(B_{\left( j+p\right)
/n})\left( |p+1|^{2H}+|p-1|^{2H}-2|p|^{2H}\right) ^{q},
\end{eqnarray*}%
which converges in $L^{1}(\Omega)$ to 
\begin{equation*}
q!\left( \sum_{k\in \mathbb{Z}}\rho _{H}(k)^q\right)
\int_{0}^{1}f(B_{s})^{2}ds\ .
\end{equation*}%
Then, it suffices to show that the terms $B_n$ and $C_n$ converge to zero in $L^2(\Omega)$.
For the term $B_{n}$ we can write, using the fact that $\sum_{l,j=0}^{n-1}%
\left| \alpha _{l,j}^{q-i}\beta _{l,j}^{i}\right| \leq Cn^{-2qH+1}$ 
\begin{eqnarray*}
E\left( \left| B_{n}\right| \right)  &\leq &Cn^{2qH-1}\sum_{i=1}^{q-1}\
\sum_{l,j=0}^{n-1}\left| \alpha _{l,j}^{q-i}\beta _{l,j}^{i}\right| \
\left\| \delta ^{q-i}\left( f^{(q-i)}(B_{l/n})\partial _{l/n}^{\otimes
(q-i)}\right) \right\| _{L^2} \\
&\leq &C\  \sum_{i=1}^{q-1}\ n^{-H(q-i)},
\end{eqnarray*}%
which converges to zero as $n$ tends to infinity. Finally, for the term $C_{n}
$ we can write%
\begin{equation*}
E\left( \left| C_{n}\right| \right) \leq \ Cn^{qH+\frac{1}{2}} \
\sup_{j}\left\| \Phi _{n}^{0,q,j}\right\| _{L^2}\leq Cn^{\frac{1}{2}-2qH+\
\left( 2H-\frac{1}{2}\right) \vee 0},
\end{equation*}%
and $\frac{1}{2}-2qH+ \left( 2H-\frac{1}{2}\right) \vee 0<0$, because if  $%
2H-\frac{1}{2}\leq 0$ this is true due to $\frac{1}{2}-2qH<0$, and if  
$2H-\frac{1}{2}\geq 0$, then we get $2H(1-q)<0$. This completes the proof
of Theorem \ref{thm1}.
\end{proof}

\subsection{Proof of the stable convergence (\ref{eq11})}

As a consequence of Theorem \ref{thm1}, we can derive the following result, which is 
nothing but (\ref{eq11}):

\begin{theorem}\label{thm-eq11}
Suppose that $f$ is a function satisfying Hypothesis \textbf{(H)}. Let $%
G_{n}$ be the sequence of random variables defined in (\ref{Gn}). Then, provided 
$H\in (\frac{1}{4q},\frac{1}{2})$, we have
$$
G_{n}-n^{-\frac{1}{2}-qH}\frac{(-1)^{q}}{2^{q}q!}\sum_{k=0}^{n-1}f^{(q)}(B_{k/n})%
\overset{\mathrm{stably}}{\longrightarrow }\,\sigma
_{H,q}\int_{0}^{1}f(B_{s})dW_{s},
$$
where $W$ is a Brownian motion independent of $B$ and $\sigma _{H,q}>0$ is
defined by (\ref{eq3a}).
\end{theorem}
\begin{proof}
\bigskip We recall first that $H_{q}\left( n^{  H}\left( \Delta B_{k/n}\right)
\right) =\frac{1}{q!}n^{qH}\delta ^{q}(\partial _{k/n}^{\otimes q})$. Then,
using (\ref{t3}) yields  
\begin{equation*}
f(B_{k/n})\delta ^{q}(\partial _{k/n}^{\otimes q})=\sum_{r=0}^{q}\binom{q}{r}%
\alpha _{k,k}^{r}\delta ^{q-r}(f^{(r)}(B_{k/n})\partial _{k/n}^{\otimes
(q-r)}),
\end{equation*}%
where $\alpha_{k,k}$ is defined in (\ref{star}).
As a consequence,%
\begin{eqnarray*}
G_{n} &=& \frac{1}{q!}n^{qH-\frac{1}{2}}\sum_{r=0}^{q}\sum_{k=0}^{n-1}%
\binom{q}{r}\alpha _{k,k}^{r}\delta ^{q-r}(f^{(r)}(B_{k/n})\partial
_{k/n}^{\otimes (q-r)}) \\
&=&\frac{1}{q!}\delta ^{q}(u_{n})+ \sum_{r=1}^{q-1} \delta
^{q-r}(v_{n}^{(r)})+R_{n},
\end{eqnarray*}%
where $u_{n}$ is defined in (\ref{v4}),%
\begin{equation*}
v_{n}^{(r)}=\frac{1}{q!}\binom{q}{r}n^{qH-\frac{1}{2}} \sum_{k=0}^{n-1}\alpha
_{k,k}^{r}\ f^{(r)}(B_{k/n})\partial _{k/n}^{\otimes (q-r)},
\end{equation*}%
and%
\begin{equation*}
R_{n}=\frac{1}{q!}n^{qH-\frac{1}{2}}\sum_{k=0}^{n-1}\alpha
_{k,k}^{q}f^{(q)}(B_{k/n}).
\end{equation*}%
The proof will be done in two steps.

\medskip

\textit{Step 1} We first show that if $H\in \left( 0,\frac{1}{2}\right) $,
and $r=1,\ldots ,q-1$,  $\delta ^{q-r}(v_{n}^{(r)})$ converges to zero in $%
L^{2}(\Omega)$ as $n$ tends to infinity. It suffices to show that $v_{n}^{(r)}$ converges
to zero in the norm of the space $\mathbb{D}^{q-r,2}(\EuFrak H^{\otimes
(q-r)})$. For $0\leq m\leq q-r$, we can write, using the
notation $\beta_{k,l}$ defined by (\ref{star}),
\begin{eqnarray*}
E\left( \left\| D^{m}v_{n}^{(r)}\right\| _{\EuFrak H^{\otimes
(q-r+m)}}^{2}\right)  &=& \left(\frac{1}{q!}\binom{q}{r}  \right)^2n^{2qH-1}  \\
&& \quad \times \sum_{k,l=0}^{n-1}\ E\left(
f^{(r+m)}(B_{k/n})f^{(r+m)}(B_{l/n})\right) \\
&&\quad \times  \alpha _{k,k}^{r}\alpha
_{l,l}^{r}\alpha _{k,l}^{m}\beta _{k,l}^{q-r} \\
&\leq &Cn^{2qH-1}n^{-2H(2r-2+m+q-r)} \\
&=&Cn^{2H-1-2Hm},
\end{eqnarray*}%
which converges to zero as $n$ tends to infinity.

\textit{Step 2 }To complete the proof it suffices to check that $$R_{n}-n^{%
-\frac{1}{2}-qH}\frac{(-1)^{q}}{2^{q}q!}\sum_{k=0}^{n-1}f^{(q)}(B_{k/n})$$
converges to zero in $L^{2}(\Omega)$ as $n$ tends to infinity. This follows from (%
\ref{d1}) and the estimates 
\begin{eqnarray*}
&&\left\| \frac{1}{q!}n^{qH-\frac{1}{2}}\sum_{k=0}^{n-1}\alpha
_{k,k}^{q}f^{(q)}(B_{k/n})-\frac{(-1)^{q}}{2^{q}q!}n^{-\frac{1}{2}%
-qH}\sum_{k=0}^{n-1}f^{(q)}(B_{k/n})\right\| _{L^2} \\
&\leq &Cn^{qH-\frac{1}{2}}\sum_{k=0}^{n-1} \left| \alpha _{k,k}^{q}-\frac{1%
}{2^{q}n^{2qH}}\right| \leq C\ n^{-qH+2H-\frac{1}{2}}.
\end{eqnarray*}%
Notice that  $-qH+2H-\frac{1}{2}<0$. The proof is now complete.
\end{proof}

\bigskip

\noindent \textbf{Acknowledgments}: We are very grateful to Giovanni
Peccati for numerous helpful comments, especially about Section \ref{sect-multiple}.
We also wish to thank an anonymous Referee for his/her careful reading.

\end{document}